\documentclass[a4paper,12pt]{article}
\usepackage{amsfonts}
\usepackage{amsmath}
\usepackage{hyperref}

\input{tcilatex}
\begin{document}

\title{Evaluating a double integral using Euler's method and Richardson
extrapolation}
\author{J. S. C. Prentice \\
%EndAName
Senior Research Officer\\
Mathsophical Ltd.\\
Johannesburg, South Africa}
\maketitle

\begin{abstract}
We transform a double integral into a second-order initial value problem,
which we solve using Euler's method and Richardson extrapolation. For an
example we consider, we achieve accuracy close to machine precision $\left(
\sim 10^{-15}\right) .$ We also use the algorithm to determine the error
curve for a Simpson cubature rule.
\end{abstract}

\section{Introduction}

Many techniques for evaluating multiple integrals have been developed (see,
for example, \cite{Stroud}$-$\cite{Prentice 2}). In this paper, we add to
this pool of knowledge: we transform a double integral into a pure
quadrature second-order initial value problem, which we solve using Euler's
method and Richardson extrapolation. We necessarily assume that the integral
has been \textquotedblleft preprocessed\textquotedblright\ (so that
singularities, discontinuities, infinite limits and so on have been dealt
with).

\section{Relevant Concepts}

Let $f\left( x,t\right) $ be a real-valued function, i.e. $f:%
%TCIMACRO{\U{211d} }%
%BeginExpansion
\mathbb{R}
%EndExpansion
^{2}\rightarrow 
%TCIMACRO{\U{211d} }%
%BeginExpansion
\mathbb{R}
%EndExpansion
,$ that is Riemann integrable in both variables. We will assume that $f$ is
suitably smooth so that all relevant derivatives in our analysis exist.

Define

\begin{equation}
C\left( x\right) \equiv \int\limits_{x_{0}}^{x}\int\limits_{t_{0}\left(
x\right) }^{t_{1}\left( x\right) }f\left( x,t\right) dtdx
\label{C(x) = double integral}
\end{equation}%
as the \textit{double integral} of interest, noting that the limits of the
inner integral may be functions of $x.$

Let $M$ denote the exact numerical value of some mathematical object, such
as an integral. Let $N\left( h\right) $ denote a numerical approximation to $%
M$ that is dependent on an adjustable parameter $h.$ Assume that 
\begin{equation}
N\left( h\right) =M+K_{1}\left( x\right) h+K_{2}\left( x\right)
h^{2}+K_{3}\left( x\right) h^{3}+K_{4}\left( x\right) h^{4}+\ldots
\label{N(h) = M + ..}
\end{equation}%
In other words, the approximation error in $N\left( h\right) $ is a power
series in $h$. \textit{Richardson extrapolation} is a process by which
values of $N\left( h\right) ,$ for differing values of $h,$ can be combined
linearly to yield approximations that are of higher order than the original
approximation $N\left( h\right) .$ Details of this procedure, relevant to
the current paper, are provided in the Appendix.

The well-known \textit{Euler's method} is a first-order Runge-Kutta method
for numerically solving an initial value problem, and it has an
approximation error of the form in (\ref{N(h) = M + ..}).

\section{Theory}

Differentiating (\ref{C(x) = double integral}) wrt $x$ gives%
\begin{align*}
C^{\prime }\left( x\right) & =\int\limits_{t_{0}\left( x\right)
}^{t_{1}\left( x\right) }f\left( x,t\right) dt \\
C^{\prime \prime }\left( x\right) & =f\left( x,t_{1}\left( x\right) \right)
t_{1}^{\prime }\left( x\right) -f\left( x,t_{0}\left( x\right) \right)
t_{0}^{\prime }\left( x\right) +\int\limits_{t_{0}\left( x\right)
}^{t_{1}\left( x\right) }\frac{\partial f\left( x,t\right) }{\partial x}dt
\end{align*}%
via the Leibniz rule.

Hence, we have%
\begin{align}
C^{\prime }& =Z  \notag \\
Z^{\prime }& =f\left( x,t_{1}\left( x\right) \right) t_{1}^{\prime }\left(
x\right) -f\left( x,t_{0}\left( x\right) \right) t_{0}^{\prime }\left(
x\right) +\int\limits_{t_{0}\left( x\right) }^{t_{1}\left( x\right) }\frac{%
\partial f\left( x,t\right) }{\partial x}dt\equiv g\left( x\right)
\label{Z'=gx}
\end{align}%
with initial values%
\begin{eqnarray*}
C\left( x_{0}\right) &=&\int\limits_{x_{0}}^{x_{0}}\int\limits_{t_{0}\left(
x\right) }^{t_{1}\left( x\right) }f\left( x,t\right) dtdx=0 \\
Z\left( x_{0}\right) &=&C^{\prime }\left( x_{0}\right)
=\int\limits_{t_{0}\left( x_{0}\right) }^{t_{1}\left( x_{0}\right) }f\left(
x_{0},t\right) dt.
\end{eqnarray*}

Of course, this is a second-order system, compactly expressed as%
\begin{equation}
\left[ 
\begin{array}{c}
C^{\prime } \\ 
Z^{\prime }%
\end{array}%
\right] =\left[ 
\begin{array}{c}
Z \\ 
g\left( x\right)%
\end{array}%
\right] .  \label{2nd order system}
\end{equation}

Euler's Method applied to this system has the form%
\begin{equation*}
\left[ 
\begin{array}{c}
C_{i+1} \\ 
Z_{i+1}%
\end{array}%
\right] =\left[ 
\begin{array}{c}
C_{i} \\ 
Z_{i}%
\end{array}%
\right] +h\left[ 
\begin{array}{c}
Z_{i} \\ 
g\left( x_{i}\right)%
\end{array}%
\right]
\end{equation*}%
where $h$ is a \textit{stepsize} - the spacing between the nodes $\left\{
x_{0},x_{1},x_{2},\ldots ,x_{i},\ldots ,x\right\} .$ We treat the upper
limit $x$ as the final node in the set. Solving this system numerically
yields the values $C\left( x_{i}\right) $ and, of course, $C\left( x\right)
. $

\section{Numerical Example}

Consider%
\begin{equation}
C\left( x\right) =\int\limits_{1}^{5}\int\limits_{x/5}^{x^{2}+1}\sin \left(
xt\right) dtdx  \label{C(x) example}
\end{equation}%
which gives, with $f\left( x,t\right) =\sin \left( xt\right) ,$%
\begin{equation}
Z^{\prime }=2x\sin \left( x^{3}+x\right) -\frac{\sin \left( \frac{x^{2}}{5}%
\right) }{5}+\int\limits_{x/5}^{x^{2}+1}t\cos tdt.
\label{Z' = ... + integral}
\end{equation}%
The one-dimensional integral in (\ref{Z' = ... + integral}) is evaluated
using composite Gaussian quadrature \cite{Prentice 1} to a precision of $%
10^{-15}$, although any suitably accurate technique can be used.

The initial values are%
\begin{align*}
C\left( 1\right) & =0 \\
Z\left( 1\right) & =\int\limits_{1/5}^{2}\sin tdt=1.396213414388384.
\end{align*}%
Clearly, we also have $x_{0}=1,x=5,t_{0}\left( x\right) =x/5$ and $%
t_{1}\left( x\right) =x^{2}+1.$

\subsection{Error control}

It is possible to choose the stepsize $h$ so as to achieve a desired level
of accuracy in the computation of $C\left( x\right) .$ Terminology and
notation used in this section is described in the Appendix, and the reader
is thereto referred.

With a moderately small value of $h$ (we used $h=0.01$ in our example), we
use Euler's method to obtain solutions of (\ref{2nd order system}), at the
nodes $\left\{ x_{0}=1,x_{1},x_{2},\ldots ,\right. $ $\left. x_{i},\ldots
,x=5\right\} ,$ for the various stepsizes required to construct $M_{4}\left(
x,h\right) $ and $M_{5}\left( x,h\right) $ via Richardson extrapolation. We
then perform the computations 
\begin{align*}
M_{4}\left( x_{i},h\right) -M_{5}\left( x_{i},h\right) & =\widetilde{K}%
_{4}\left( x_{i}\right) h^{4}-O\left( h^{5}\right) \\
& \approx \widetilde{K}_{4}\left( x_{i}\right) h^{4} \\
\Rightarrow \widetilde{K}_{4}\left( x_{i}\right) & =\frac{M_{4}\left(
x_{i},h\right) -M_{5}\left( x_{i},h\right) }{h^{4}}
\end{align*}%
at each of the nodes $\left\{ x_{0}=1,x_{1},x_{2},\ldots ,x_{i},\ldots
,x=5\right\} $. Hence, we have the error coefficient $\widetilde{K}_{4}$ at
each node. We note that, due to the linear combination necessary to
construct $M_{4}\left( x,h\right) ,$ $\widetilde{K}_{4}$ is not necessarily
equal to $K_{4}.$ In fact, we have $K_{4}=64\widetilde{K}_{4}.$

The sequence of calculations below then allows a new stepsize to be found,
consistent with a user-defined tolerance $\varepsilon .$ 
\begin{align*}
h^{\ast }& =\left( \frac{\varepsilon }{\max\limits_{x_{i}}\left\vert 
\widetilde{K}_{4}\left( x_{i}\right) \right\vert }\right) ^{\frac{1}{4}} \\
n& =\left\lceil \frac{x-x_{0}}{h^{\ast }}\right\rceil \\
h& =\frac{x-x_{0}}{n}.
\end{align*}%
The Euler/Richardson algorithm is then repeated using this new stepsize,
ultimately leading to a new $M_{4}\left( x_{i},h\right) $ and, in
particular, $M_{4}\left( 5,h\right) \approx C\left( 5\right) .$

\subsection{Results}

Applying the algorithm to (\ref{C(x) example}) gives the stepsizes needed
for various tolerances, shown in Table 1.

\medskip

\begin{center}
\begin{tabular}{|l||l|l|l|l|l|l|}
\multicolumn{7}{l}{Table 1: Tolerances and corresponding stepsizes} \\ \hline
$\varepsilon $ & $10^{-14}$ & $10^{-12}$ & $10^{-10}$ & $10^{-8}$ & $10^{-6}$
& $10^{-4}$ \\ \hline
$h$ & $2\times 10^{-4}$ & $6.3\times 10^{-4}$ & $2\times 10^{-3}$ & $%
6.3\times 10^{-3}$ & $2\times 10^{-2}$ & $6.3\times 10^{-2}$ \\ \hline
\end{tabular}
\end{center}

\medskip

In Figure 1 we show $C\left( x\right) $ obtained from $M_{4}\left(
x_{i},h\right) $ with $h=6.3\times 10^{-4}.$ In Figure 2, we show the
associated coefficients $K_{1}\left( x\right) ,K_{2}\left( x\right)
,K_{3}\left( x\right) ,K_{4}\left( x\right) $ and $K_{5}\left( x\right) .$
In Figure 3, we show error curves for $C\left( x\right) $ obtained from both 
$M_{2}\left( x_{i},h\right) $ and $M_{4}\left( x_{i},h\right) ,$ again with $%
h=6.3\times 10^{-4}.$ Clearly, $M_{4}\left( x_{i},h\right) $ appears to
achieve a tolerance of $\sim 10^{-12}$ over the entire interval.

Of course, the original objective was to determine $C\left( 5\right) ,$ and
so%
\begin{equation*}
C\left( 5\right) \approx M_{4}\left( 5,6.3\times 10^{-4}\right)
=0.630635228375177
\end{equation*}%
which is within $8.3\times 10^{-13}$ of the exact value.

\section{Possible Applications}

Apart from the obvious task of determining $C\left( x\right) ,$ we can also
use this algorithm to find the error term for a given cubature rule $Q\left(
x\right) .$ If we include $Q\left( x\right) $ in (\ref{C(x) = double
integral}), we have 
\begin{equation*}
\int\limits_{x_{0}}^{x}\int\limits_{t_{0}\left( x\right) }^{t_{1}\left(
x\right) }f\left( x,t\right) dtdx=Q\left( x\right) +C\left( x\right) 
\end{equation*}%
which yields the system%
\begin{align*}
C^{\prime }& =Z \\
Z^{\prime }& =f\left( x,t_{1}\left( x\right) \right) t_{1}^{\prime }\left(
x\right) -f\left( x,t_{0}\left( x\right) \right) t_{0}^{\prime }\left(
x\right)  \\
& +\int\limits_{t_{0}\left( x\right) }^{t_{1}\left( x\right) }\frac{\partial
f\left( x,t\right) }{\partial x}dt-Q^{\prime \prime }\left( x\right)  \\
& \equiv g\left( x\right) ,
\end{align*}%
with initial values%
\begin{align*}
C\left( x_{0}\right) & =\int\limits_{x_{0}}^{x_{0}}\int\limits_{t_{0}\left(
x\right) }^{t_{1}\left( x\right) }f\left( x,t\right) dtdx-Q\left(
x_{0}\right) =-Q\left( x_{0}\right)  \\
Z\left( x_{0}\right) & =C^{\prime }\left( x_{0}\right)
=\int\limits_{t_{0}\left( x_{0}\right) }^{t_{1}\left( x_{0}\right) }f\left(
x_{0},t\right) dt-Q^{\prime }\left( x_{0}\right) .
\end{align*}%
In this case, $C\left( x\right) $ now plays the role of an \textit{error term%
} or a \textit{correction term} (earlier, we effectively had $Q\left(
x\right) =0).$

For example, applying Simpson's Rule to our example yields the cubature
expression $Q\left( x\right) $ shown in the Appendix. Then, applying our
algorithm with $h=4\times 10^{-3},$ we find the curves $Q\left( x\right) $
and $C\left( x\right) $ shown in Figure 4. Comparing $Q\left( x\right) $ to
the curve in Figure 1, we see that the cubature rule is very inaccurate.
However, when the correction term $C\left( x\right) $ is added to $Q\left(
x\right) ,$ the result differs from the true value by $\sim 2\times 10^{-6}$
(for $M_{4})$ and $\sim 5\times 10^{-9}$ (for $M_{5})$ at $x=5,$ and these
errors can, of course, be made even smaller by using a smaller stepsize.

\section{Conclusion}

We have described a technique by which the task of determining a double
integral can be presented as the task of solving a second-order initial
value problem. We use Euler's method for this purpose, combined with
Richardson extrapolation to yield very accurate results. Indeed, we have
achieved accuracy close to machine precision for the numerical example
considered here. Moreover, our algorithm lends itself to error control - the
stepsize $h$ needed for a desired accuracy can be estimated. Additionally,
we have shown how the algorithm can be adapted to yield the error curve for
a cubature rule, if such a rule is imposed on the problem.

Future work could address the matter of relative error control,
higher-dimensional cubature and the possibility of developing custom
Runge-Kutta methods of high order that are able to incorporate the integral
term that appears in (\ref{Z'=gx}). If the interval of integration is
relatively large, it may be wise to divide said interval into subintervals
and to apply the algorithm on each subinterval, adding the results in an
appropriate manner. Whether or not this is necessary or advantageous needs
to investigated. The general case where the limits in the outer integral are
arbitrary functions of $x$ should also be considered (we make some comments
in this regard in the Appendix).

Lastly, we note that the physical time taken on our platform \cite{platform}
to compute $C\left( x\right) $ - for $\sim 6100$ nodes - was roughly $30$s.
To compute the same curve using the symbolic \texttt{int(f,a,x)} function in
MatLab took as long as an hour. Future research might study this aspect of
efficiency in more detail.
\begin{verbatim}
 
\end{verbatim}

\section{\protect\medskip Appendix}

\subsection{\protect\medskip Richardson extrapolation}

We define%
\begin{align*}
N_{0}\left( x,h\right) \equiv N\left( x,h\right) =N\left( x,\frac{h}{2^{0}}%
\right) =& M\left( x\right) +K_{1}\left( x\right) h+K_{2}\left( x\right)
h^{2} \\
& +K_{3}\left( x\right) h^{3}+K_{4}\left( x\right) h^{4}+K_{5}\left(
x\right) h^{5}+\ldots \\
N_{1}\left( x,h\right) \equiv N\left( x,\frac{h}{2}\right) =N\left( x,\frac{h%
}{2^{1}}\right) =& M\left( x\right) +K_{1}\left( x\right) \frac{h}{2}%
+K_{2}\left( x\right) \frac{h^{2}}{4} \\
& +K_{3}\left( x\right) \frac{h^{3}}{8}+K_{4}\left( x\right) \frac{h^{4}}{16}%
+K_{5}\left( x\right) \frac{h^{5}}{32}+\ldots .
\end{align*}%
and similarly for 
\begin{align*}
N_{2}\left( x,h\right) & \equiv N\left( x,\frac{h}{2^{2}}\right) =N\left( x,%
\frac{h}{4}\right) \\
N_{3}\left( x,h\right) & \equiv N\left( x,\frac{h}{2^{3}}\right) =N\left( x,%
\frac{h}{8}\right) \\
N_{4}\left( x,h\right) & \equiv N\left( x,\frac{h}{2^{4}}\right) =N\left( x,%
\frac{h}{16}\right) \\
N_{5}\left( x,h\right) & \equiv N\left( x,\frac{h}{2^{5}}\right) =N\left( x,%
\frac{h}{32}\right) .
\end{align*}

For example, to construct a fifth-order method, we seek $\alpha _{0},\ldots
,\alpha _{4}$ such that%
\begin{equation*}
M_{5}\left( x,h\right) \equiv \sum\limits_{k=0}^{4}\alpha _{k}N_{k}\left(
x,h\right) =M\left( x\right) +O\left( h^{5}\right) .
\end{equation*}

It can be shown that we must solve the linear system\renewcommand{%
\arraystretch}{1.7}%
\begin{equation*}
\left[ 
\begin{array}{ccccc}
1 & 1 & 1 & 1 & 1 \\ 
1 & \frac{1}{2} & \frac{1}{4} & \frac{1}{8} & \frac{1}{16} \\ 
1 & \frac{1}{4} & \frac{1}{16} & \frac{1}{64} & \frac{1}{256} \\ 
1 & \frac{1}{8} & \frac{1}{64} & \frac{1}{512} & \frac{1}{4096} \\ 
1 & \frac{1}{16} & \frac{1}{256} & \frac{1}{4096} & \frac{1}{65536}%
\end{array}%
\right] \left[ 
\begin{array}{c}
\alpha _{0} \\ 
\alpha _{1} \\ 
\alpha _{2} \\ 
\alpha _{3} \\ 
\alpha _{4}%
\end{array}%
\right] =\left[ 
\begin{array}{c}
1 \\ 
0 \\ 
0 \\ 
0 \\ 
0%
\end{array}%
\right] ,
\end{equation*}%
\renewcommand{\arraystretch}{1}which gives\renewcommand{\arraystretch}{1.7}%
\begin{equation*}
\left[ 
\begin{array}{c}
\alpha _{0} \\ 
\alpha _{1} \\ 
\alpha _{2} \\ 
\alpha _{3} \\ 
\alpha _{4}%
\end{array}%
\right] =\left[ 
\begin{array}{r}
\frac{1}{315} \\ 
-\frac{2}{21} \\ 
\frac{8}{9} \\ 
-\frac{64}{21} \\ 
\frac{1024}{315}%
\end{array}%
\right] .
\end{equation*}%
\renewcommand{\arraystretch}{1}

In general, we have, for an $m$th order method 
\begin{equation*}
M_{m}\left( x,h\right) \equiv \sum\limits_{k=0}^{m-1}\alpha _{k}N_{k}\left(
x,h\right) =M\left( x\right) +O\left( h^{m}\right) ,
\end{equation*}%
\begin{equation*}
\left[ 
\begin{array}{ccccc}
A_{1,1} & \cdots & \cdots & \cdots & A_{1,m} \\ 
\vdots & \ddots &  &  & \vdots \\ 
\vdots &  & A_{i,j} &  & \vdots \\ 
\vdots &  &  & \ddots & \vdots \\ 
A_{m,1} & \cdots & \cdots & \cdots & A_{m,m}%
\end{array}%
\right] \left[ 
\begin{array}{c}
\alpha _{0} \\ 
\vdots \\ 
\vdots \\ 
\vdots \\ 
\alpha _{m-1}%
\end{array}%
\right] =\left[ 
\begin{array}{c}
1 \\ 
0 \\ 
\vdots \\ 
\vdots \\ 
0%
\end{array}%
\right] \text{ \ where \ }A_{i,j}=\left( \frac{1}{2^{i-1}}\right) ^{j-1}.
\end{equation*}

Using this, we can find coefficients for 2nd-, 3rd- and 4th-order methods%
\renewcommand{\arraystretch}{1.7}%
\begin{equation*}
\left[ 
\begin{array}{c}
\alpha _{0} \\ 
\alpha _{1}%
\end{array}%
\right] =\left[ 
\begin{array}{r}
-1 \\ 
2%
\end{array}%
\right] ,\text{ \ \ \ \ \ \ }\left[ 
\begin{array}{c}
\alpha _{0} \\ 
\alpha _{1} \\ 
\alpha _{2}%
\end{array}%
\right] =\left[ 
\begin{array}{r}
\frac{1}{3} \\ 
-2 \\ 
\frac{8}{3}%
\end{array}%
\right] ,\text{ \ \ \ \ \ \ }\left[ 
\begin{array}{c}
\alpha _{0} \\ 
\alpha _{1} \\ 
\alpha _{2} \\ 
\alpha _{3}%
\end{array}%
\right] =\left[ 
\begin{array}{r}
-\frac{1}{21} \\ 
\frac{2}{3} \\ 
-\frac{8}{3} \\ 
\frac{64}{21}%
\end{array}%
\right]
\end{equation*}%
\renewcommand{\arraystretch}{1}and a 6th-order method\renewcommand{%
\arraystretch}{1.7}%
\begin{equation*}
\left[ 
\begin{array}{c}
\alpha _{0} \\ 
\alpha _{1} \\ 
\alpha _{2} \\ 
\alpha _{3} \\ 
\alpha _{4} \\ 
\alpha _{5}%
\end{array}%
\right] =\left[ 
\begin{array}{r}
-\frac{1}{9765} \\ 
\frac{2}{315} \\ 
-\frac{8}{63} \\ 
\frac{64}{63} \\ 
-\frac{1024}{315} \\ 
\frac{32768}{9765}%
\end{array}%
\right] .
\end{equation*}%
\renewcommand{\arraystretch}{1}

\subsection{Simpson cubature}

$\medskip $A cubature expression, in terms of $x,$ obtained by applying
Simpson's Rule to (\ref{C(x) example}), is

\medskip

\begin{center}
\begin{tabular}{ll}
$Q(x)$ $\ =$ & $\frac{\left( x-1\right) }{360}\left( 18\sin \left( 2\right)
+18\,\sin \left( \frac{1}{5}\right) +72\,\sin \left( \frac{11}{10}\right)
\right) +$ \\ 
& $\frac{\left( x-1\right) }{36}\left( x^{2}-\frac{x}{5}+1\right) \,\left( 
\begin{array}{c}
\sin \left( x\,\left( x^{2}+1\right) \right) +\sin \left( \frac{x^{2}}{5}%
\right) + \\ 
4\,\sin \left( x\,\left( \frac{x^{2}}{2}+\frac{x}{10}+\frac{1}{2}\right)
\right)%
\end{array}%
\right) +$ \\ 
& $\frac{\left( x-1\right) }{36}\left( {\left( \frac{x}{2}+\frac{1}{2}%
\right) }^{2}-\frac{x}{10}+\frac{9}{10}\right) \,\left( 
\begin{array}{c}
16\,\sin \left( \left( \frac{x}{2}+\frac{1}{2}\right) \,\left( \frac{x}{20}+%
\frac{{\left( \frac{x}{2}+\frac{1}{2}\right) }^{2}}{2}+\frac{11}{20}\right)
\right) + \\ 
4\,\sin \left( \left( \frac{x}{2}+\frac{1}{2}\right) \,\left( \frac{x}{10}+%
\frac{1}{10}\right) \right) + \\ 
4\,\sin \left( \left( \frac{x}{2}+\frac{1}{2}\right) \,\left( {\left( \frac{x%
}{2}+\frac{1}{2}\right) }^{2}+1\right) \right)%
\end{array}%
\right) .$%
\end{tabular}

\medskip
\end{center}

This expression was obtained using symbolic software \cite{platform}, and
includes the effect of the variable limits in the inner integral.

\subsection{General limits}

Consider the more general case%
\begin{equation*}
C\left( x\right) \equiv \int\limits_{a\left( x\right) }^{b\left( x\right)
}\int\limits_{t_{0}\left( x\right) }^{t_{1}\left( x\right) }f\left(
x,t\right) dtdx\equiv \int\limits_{a\left( x\right) }^{b\left( x\right)
}I\left( x\right) dx
\end{equation*}%
where the limits in the outer integral are both functions of $x,$ not
necessarily linear, and we have implicitly defined $I\left( x\right) .$

We find, dropping the argument of $a\left( x\right) $ and $b\left( x\right) $
for notational convenience,%
\begin{align}
C^{\prime }\left( x\right) =\text{ }& I\left( b\left( x\right) \right)
b^{\prime }\left( x\right) -I\left( a\left( x\right) \right) a^{\prime
}\left( x\right)  \notag \\
C^{\prime \prime }\left( x\right) =\text{ }& I\left( b\left( x\right)
\right) b^{\prime \prime }\left( x\right) -I\left( a\left( x\right) \right)
a^{\prime \prime }\left( x\right)  \label{C'' = I(b)B''...} \\
& +\frac{dI\left( b\right) }{db}b^{\prime }\left( x\right) b^{\prime }\left(
x\right) -\frac{dI\left( a\right) }{da}a^{\prime }\left( x\right) a^{\prime
}\left( x\right) .  \notag
\end{align}

Now,%
\begin{align*}
\frac{dI\left( b\right) }{db}& =\frac{d}{db}\left( \int\limits_{t_{0}\left(
b\right) }^{t_{1}\left( b\right) }f\left( b,t\right) dt\right) \\
& =f\left( b,t_{1}\left( b\right) \right) \frac{dt_{1}\left( b\right) }{db}%
-f\left( b,t_{0}\left( b\right) \right) \frac{dt_{0}\left( b\right) }{db}%
+\int\limits_{t_{0}\left( b\right) }^{t_{1}\left( b\right) }\frac{\partial
f\left( b,t\right) }{\partial b}dt
\end{align*}%
and%
\begin{align*}
\frac{dI\left( a\right) }{da}& =\frac{d}{da}\left( \int\limits_{t_{0}\left(
a\right) }^{t_{1}\left( a\right) }f\left( a,t\right) dt\right) \\
& =f\left( a,t_{1}\left( a\right) \right) \frac{dt_{1}\left( a\right) }{da}%
-f\left( a,t_{0}\left( a\right) \right) \frac{dt_{0}\left( a\right) }{da}%
+\int\limits_{t_{0}\left( a\right) }^{t_{1}\left( a\right) }\frac{\partial
f\left( a,t\right) }{\partial a}dt.
\end{align*}

There are two approaches we may take to solving (\ref{C'' = I(b)B''...}).
The first is the case where a specific value of $x$ is given, say $x=p$. We
then simply use $p$ to determine the numerical values for the outer limits,
giving%
\begin{equation*}
C\left( x\right) \equiv \int\limits_{a\left( p\right) }^{b\left( p\right)
}\int\limits_{t_{0}\left( x\right) }^{t_{1}\left( x\right) }f\left(
x,t\right) dtdx.
\end{equation*}%
This has the same form as (\ref{C(x) example}), and may be handled in the
same way. The second case is when we desire $C\left( x\right) $ over an
interval, say $\left[ x_{0},x_{n}\right] .$ In this case, we solve the system%
\begin{align*}
C^{\prime }=\text{ }& Z \\
Z^{\prime }=\text{ }& I\left( b\left( x\right) \right) b^{\prime \prime
}\left( x\right) -I\left( a\left( x\right) \right) a^{\prime \prime }\left(
x\right)  \\
& +\frac{dI\left( b\right) }{db}b^{\prime }\left( x\right) b^{\prime }\left(
x\right) -\frac{dI\left( a\right) }{da}a^{\prime }\left( x\right) a^{\prime
}\left( x\right) 
\end{align*}%
with initial values%
\begin{align*}
C\left( x_{0}\right) & =\int\limits_{a\left( x_{0}\right) }^{b\left(
x_{0}\right) }I\left( x_{0}\right) dx=\left( b\left( x_{0}\right) -a\left(
x_{0}\right) \right) I\left( x_{0}\right)  \\
Z\left( x_{0}\right) & =I\left( b\left( x_{0}\right) \right) b^{\prime
}\left( x_{0}\right) -I\left( a\left( x_{0}\right) \right) a^{\prime }\left(
x_{0}\right) .
\end{align*}%
We will not investigate this line of research any further here, preferring
to consider it in future work.

Nevertheless, for the sake of completeness: in our earlier example we had%
\begin{align*}
a\left( x\right) & =x_{0}=const. \\
b\left( x\right) & =x
\end{align*}%
so that%
\begin{align*}
a^{\prime }& =a^{\prime \prime }=0 \\
b^{\prime }& =1,\text{ }b^{\prime \prime }=0 \\
t_{1}\left( b\right) & =t_{1}\left( x\right) ,\text{ }t_{0}\left( b\right)
=t_{0}\left( x\right) \\
t_{1}\left( a\right) & =t_{1}\left( x\right) ,\text{ }t_{0}\left( a\right)
=t_{0}\left( x\right) .
\end{align*}%
Substituting these into (\ref{C'' = I(b)B''...}) gives%
\begin{align*}
C^{\prime }\left( x\right) & =I\left( b\left( x\right) \right) b^{\prime
}\left( x\right) -I\left( a\left( x\right) \right) a^{\prime }\left(
x\right) =I\left( x\right) =\int\limits_{t_{0}\left( x\right) }^{t_{1}\left(
x\right) }f\left( x,t\right) dt \\
C^{\prime \prime }\left( x\right) & =0-0+\frac{dI\left( x\right) }{dx}\left(
1\right) \left( 1\right) -0 \\
& =f\left( x,t_{1}\left( x\right) \right) \frac{dt_{1}\left( x\right) }{dx}%
-f\left( x,t_{0}\left( x\right) \right) \frac{dt_{0}\left( x\right) }{dx}%
+\int\limits_{t_{0}\left( x\right) }^{t_{1}\left( x\right) }\frac{\partial
f\left( x,t\right) }{\partial x}dt,
\end{align*}%
as expected.

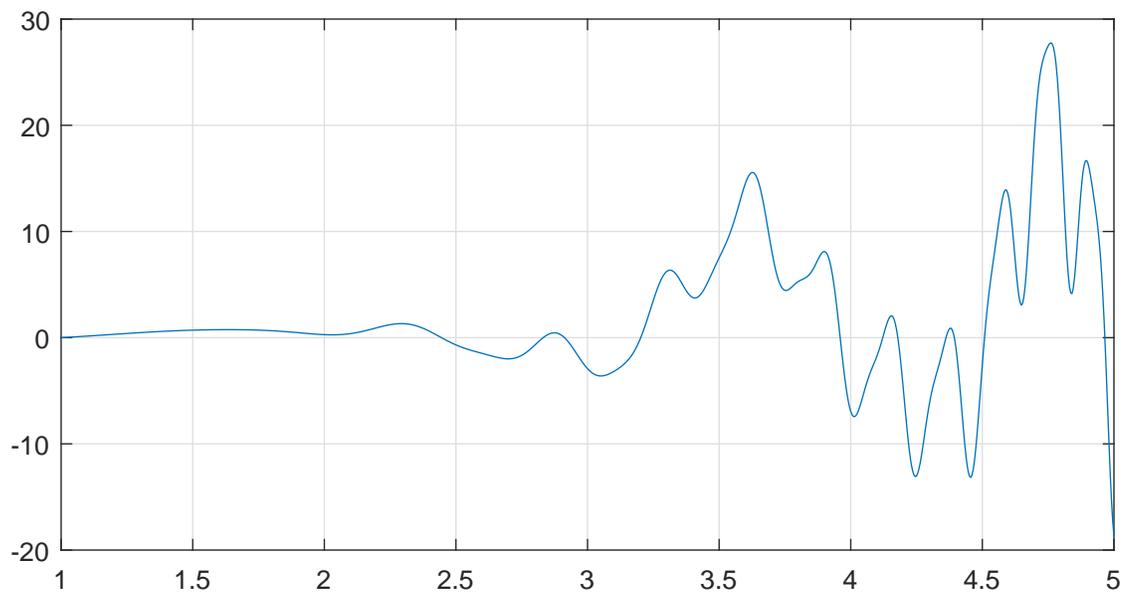
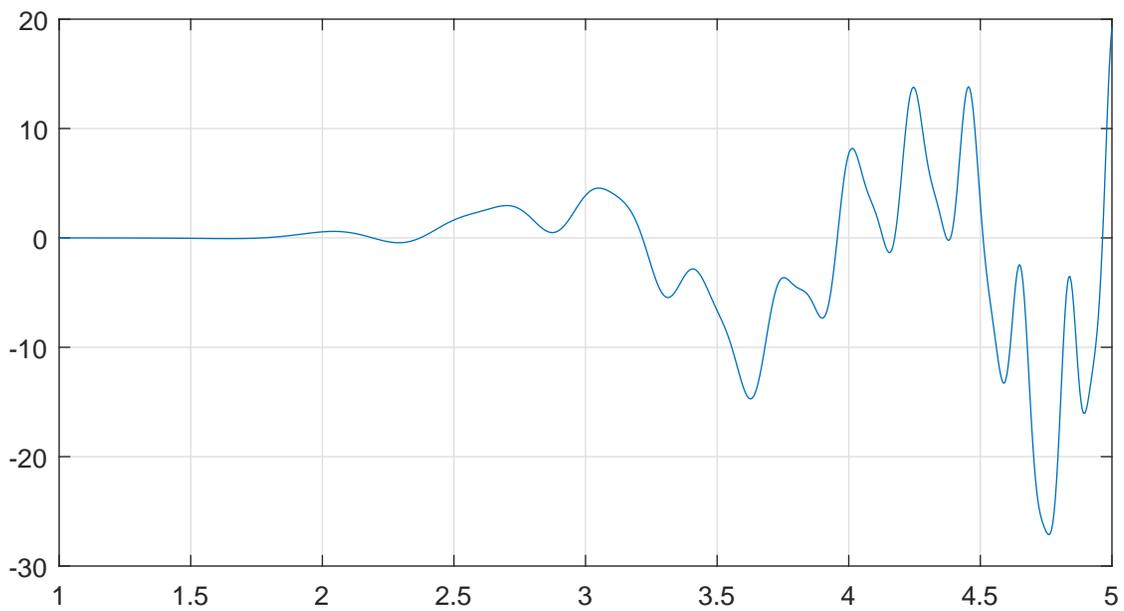
\end{document}